\input amstex
\documentstyle{amsppt}

\nologo
\parskip=6pt
\magnification=1200

\input epsf

\topmatter
\title KNOT ENERGIES AND KNOT INVARIANTS
\endtitle
\author Xiao-Song Lin
\endauthor
\affil Department of Mathematics \\ University of California 
\\ Riverside, CA 92521
\\{\it xl\@math.ucr.edu}
\endaffil
\endtopmatter 

\document

\rightline{$\ssize\text{To record what has happened, ancient people tie
knots.}$}
\rightline{$\ssize\text{--- {\it I ching}, the Chinese classic of 1027--771 
B.C.}$}
\bigskip

Knots are fascinating objects. When fastening a rope, the distinction between
a knot and a \lq\lq slip-knot" (one that can be
undone by pulling) must be recognized very
early in human history. We even developed a subconscious about knots: When we 
are puzzled or troubled, we have a feeling of being knotted somewhere. The
mathematical study of knots started much later though. It was inspired in the
middle of the nineteenth century by the vortex theory of fluid dynamics (see
[11] for a vivid description of this history). The development of
modern topology in the first half of the twentieth century provided a solid
background for a mathematical theory of knots. Yet we only began to see the
full scope of knot theory in the last decade, starting with the discovery of 
the
Jones polynomial in 1984 (see [13] for a survey of  
the history of knot theory up to Jones' discovery). 
In 1989, Witten generalized the Jones polynomial using his
Chern-Simons path integral. Finally, in 1990-92, the development of knot theory
culminated in the theory of Vassiliev knot invariants, which provides probably
the most general framework for the study of the combinatorics of knots.
Through
the study of Vassiliev knot invariants, we see that although 
the abundance of knots in varieties is distinctively visible, this
abundance does not come from any randomness. The
combinatorics of knots embraces almost all fundamental symmetries of 
mathematics
and physics that we know. Such a pervasive nature is not common among
topological and geometric objects mathematicians are in favor of. For the 
reader's convenience, we have collected several excellent expository papers
on these developments in the references (see [1, 2, 4, 8, 14]).

Geometers are restless in their effort of searching for geometric objects with
``maximal homogeneity". Here,
of course, the measurement of homogeneity are different in different
situations. Actually, it is the key point to recognize in a given geometric
setting what should be the measurement of homogeneity. Thus, in classical
Riemannian geometry, we know that various curvatures are the key measurement of
homogeneity; we measure length or area for immersions of circles and surfaces 
into a Riemannian manifold and developed the theories of geodesics and minimal
surfaces; in gauge theory, we study connections minimizing the Yang-Mills
functional; and we look for pseudo-holomorphic
curves in symplectic geometry; etc.. And there is always the moduli problem if
geometric objects with maximal homogeneity are not unique. So, we may also 
ask for smooth imbeddings
$S^1\rightarrow \Bbb R^3$, which we will refer to as {\it geometric knots}, 
with
the ``most perfect" shape among all geometric knots isotopic to each other.
This geometric side of knot theory is much less mature than the combinatorial
side of knot theory. It seems to be not completely clear yet as of what should
be the most fundamental measurement for the homogeneity of a geometric knot
within its isotopy class. One of the purposes of this article is then to argue
that such measurements of homogeneity satisfying the criteria set forth in the
foundational paper of Freedman, He and Wang [6] may not be unique. There 
seems to be
a spectrum of M\"obius knot energies related with
functionals on geometric knots appeared in integral formulae of
Vassiliev knot invariants coming from perturbative expansion of Witten's
Chern-Simons path integral (they are called {\it Gauss functionals}). 
Of course, unless we could understand the
dynamical behavior of geometric knots with respect to these M\"obius energies, 
their nature remains to be a mystery. 

Classically, functionals on loop spaces people have studied include the
length functional and holonomy functionals. Functionals defined only on 
embedded curves have caught people's attention lately, to a large part
due to the recent advance in knot theory. The elementary
discussion on
those Gauss functionals on geometric knots in this article reminds us
classical integral geometry where different 
measurements on the same geometric object are shown to be related. 
Hopefully, this will motivate further interesting in geometric knot 
theory. 

We would like to thank Colin Adams for the invitation of writing an article for
this special issue of Chaos, Solitons and Fractals. 
We have talked about the topics of this 
article in geometry/topology seminars at UCD, UCR and UCSD. It is a pleasure
to thank the participants of these seminars, in particular, 
Professors Mike Freedman and Oleg Viro, for their interests and 
comments. It would be impossible to have the work presented here done 
without continuing discussions with Zhenghan Wang. We appreciate very much the
help and stimulation we get from him.      

\bigskip
\noindent{\bf \S1. M\"obius energies}
\medskip

We define a {\it geometric knot} to be a smooth embedding
$\gamma:S^1\rightarrow\Bbb R^3$. Here, the oriented circle $S^1$ comes with no
particular parameterization. 

Two geometric knots $\gamma_1$ and $\gamma_2$ are {\it equivalent} (or {\it
isotopic}) if there is an orientation preserving diffeomorphism $\rho:\Bbb
R^3\rightarrow\Bbb R^3$ such that $\gamma_2=\rho\gamma_1$. A {\it knot} is
simply an equivalence class of geometric knots. See Figure 1. So, a geometric 
knot is {\it unknotted} if it can be deformed to a round circle without passing
through itself.

The {\it M\"obius energy} $E(\gamma)$ of a geometric knot $\gamma$, defined
in [6], is given as follows: Suppose $S^1$ is parameterized by $u$ and 
$du$ is
the positive volume form of $S^1$ coming from that 
parameterization. For $u,v\in
S^1$, $u\neq v$, denote by $D(\gamma(v),\gamma(u))$ be the minimum of the
lengths of sub-arcs of $\gamma$ from $\gamma(u)$ to $\gamma(v)$.
Then
$$E(\gamma)=\int_{u\in S^1}|\dot\gamma(u)|\,du\left(\int_{v\in S^1,
v\neq u}\left
\{\frac1{|\gamma(v)-\gamma(u)|^2}-
\frac1{D(\gamma(v),\gamma(u))^2}\right\}|\dot\gamma(v)|\,dv\right).$$
The integral is independent of the parameterization of $S^1$. It is therefore a
positive functional on geometric knots. It is also independent of
the orientation of $S^1$.

\bigskip
\centerline{\epsfxsize=2.5in \epsfbox{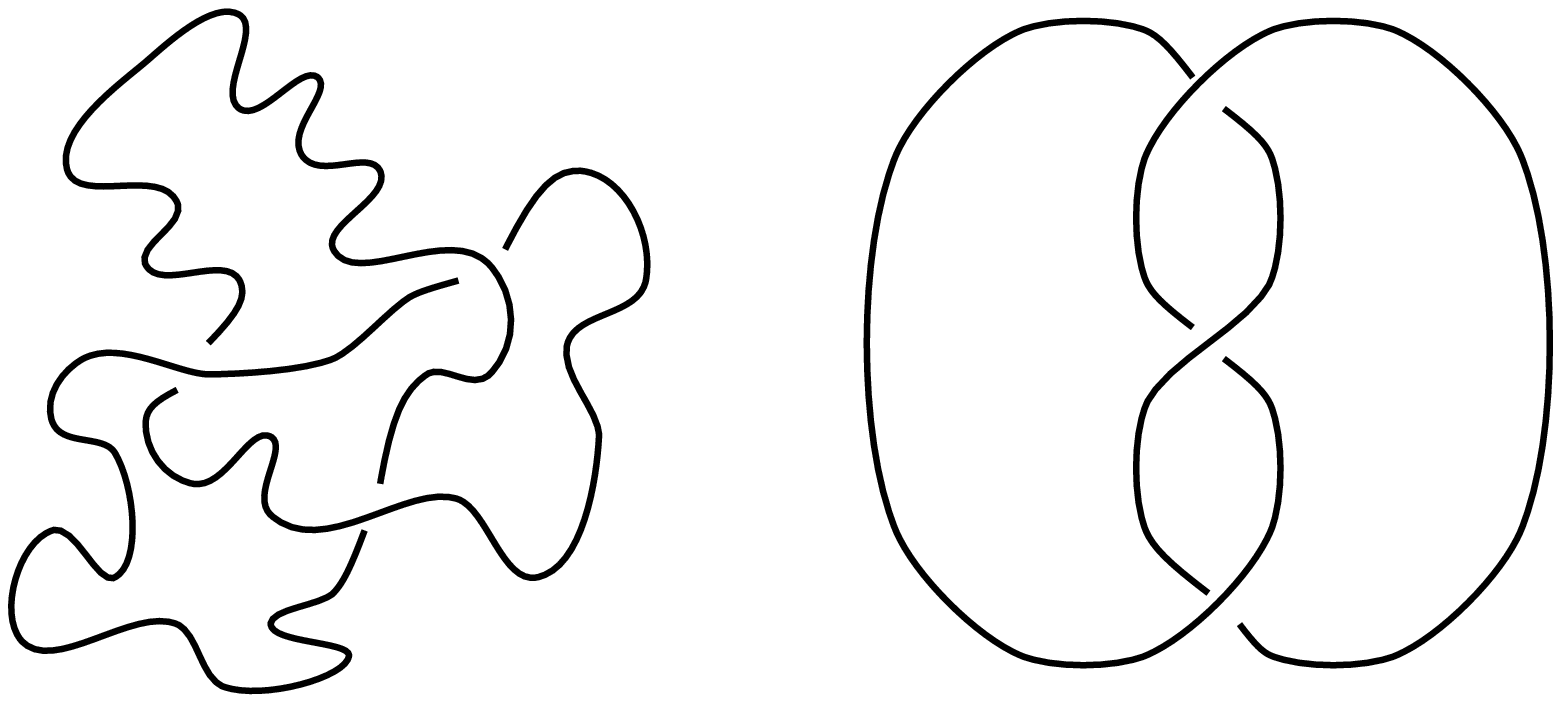}}
\smallskip
\centerline{$\ssize\text{Figure 1. Two equivalent (isotopic) geometric knot.}$}
\medskip

The idea is that as a geometric knot $\gamma$ deforms and tends to acquire a
double point,
$E(\gamma)$ will blow up and thus constraints the deformation of $\gamma$ 
within its isotopy class. Ideally, every geometric knot 
would deform to a unique energy-minimizing geometric knot within 
its isotopy class via
the gradient flow of $E$ (the direction where $E$ is decreasing). Although 
this is not true in general, 
we will see below that geometric knots with smaller $E$ do
look more ``homogeneous''. See [9] for a table of many 
energy-minimizing 
geometric knots obtained by computer simulation. 

A fundamental property of the energy functional $E$ discovered in [6] is its
invariance under M\"obius transformation of $\Bbb R^3\cup\{\infty\}$.
M\"obius transformations of $\Bbb R^3\cup\{\infty\}$ are the 10-dimensional Lie
group of angle-preserving diffeomorphisms of $\Bbb R^3\cup\{\infty\}$ generated
by inversion in 2-spheres. The M\"obius invariance of $E$ says that if $T$
is a M\"obius transformation and $T\gamma\subset\Bbb R^3$, then 
$E(T\gamma)=E(\gamma)$. If $T\gamma\subset\Bbb R^3\cup\{\infty\}$, one
may generalize the energy functional $E$ to infinite curves in $\Bbb R^3$ and
then we have
$E(T\gamma\cap\Bbb R^3)=E(\gamma)-4$. See Theorem 2.1 of [6]. 

The M\"obius invariance of $E$ led to its reformulation by P. Doyle and
O. Schramm in terms of some conformally invariant data.  

Let $u,v\in S^1$, $u\neq v$. In $\Bbb R^3$, there are exactly two round, 
oriented circles or oriented straight lines $U_v$ and $V_u$ passing through
$\gamma(u)$ and $\gamma(v)$ and tangent to 
$\dot\gamma(u)$ at $\gamma(u)$ and $\dot\gamma(v)$ at
$\gamma(v)$ respectively. Let $\alpha=\alpha(u,v)$ be the angle between $U_u$ 
and $V_u$, $0\leq\alpha\leq\pi$. Define 
$$
E_{\cos}(\gamma)=\int_{u\in S^1}|\dot\gamma(u)|\,du\left(
\int_{v\in S^1,v\neq
u}\frac{1-\cos\alpha}{|\gamma(v)-\gamma(u)|^2}|\dot\gamma(v)|\,dv\right).$$
The M\"obius invariance of $E_{\cos}$ is more or less transparent.

The functional $E_{\cos}$ may be interpreted in terms of 
``excess lengths''.
Fix $u\in S^1$ and we may assume that $\gamma(u)=0$ and $\dot\gamma(u)$
is horizontal. We apply the M\"obius inversion
$x\mapsto\dsize\frac{x}{|x|^2}$ about the unit 2-sphere centered at the origin
to $\gamma$. Then $\gamma$ becomes an asymptotically horizontal infinite 
curve $\gamma_{\infty}$ in $\Bbb R^3$. The M\"obius inversion sends each
$U_v$ to a horizontal straight line $L_v$ and they are all parallel to 
each other for different $v$'s. The angle $\alpha$ between $U_v$ and $V_u$
is the same as the angle between $\dot\gamma_\infty(v)$ and $L_v$. Although
both $\gamma_\infty$ and a horizontal straight line have infinite lengths, 
these two infinities are comparable in the sense that their difference can be
made finite. This excess length can be computed in the following way. 

Let $s$ be the arc length parameter of $\gamma_\infty$. 
Then, from $s_1$ to $s_2$, the horizontal distance one travels along 
$\gamma_\infty$ is
$$\int_{s_1}^{s_2}\cos\alpha\,ds.$$
Therefore the {\it horizontal excess length} of $\gamma_\infty$ is
$$\int_{\gamma_\infty}(1-\cos\alpha)\,ds<\infty.$$
But 
$$ds=|\dot\gamma_\infty(v)|\,dv\qquad\text{and}\qquad
\gamma_\infty(v)=\frac{\gamma(v)}{|\gamma(v)|^2}.$$
Moreover, simple vector calculus shows
$$|\dot\gamma_\infty(v)|=\frac{|\dot\gamma(v)|}{|\gamma(v)|^2}.$$
Thus, the first integration in $E_{\cos}$ is exactly the 
horizontal excess length
of $\gamma_{\infty}$. From this interpretation, it is quite clear that a
geometric knot which has a ``highly oscillated segment'' can not be 
a $E_{\cos}$-minimizing geometric knot.

\proclaim{\bf Proposition} {\rm (Doyle-Schramm)} $E=E_{\cos}+4$.
\endproclaim

There is one more M\"obius energy paired with $E_{\cos}$:
$$E_{\sin}(\gamma)=\int_{u\in S^1}|\dot\gamma(u)|\,du\left(
\int_{v\in S^1,v\neq
u}\frac{\sin\alpha}{|\gamma(v)-\gamma(u)|^2}|\dot\gamma(v)|\,dv\right).$$

The energy $E_{\sin}$ lacks the smoothness of $E_{\cos}$. To see this, 
let $\alpha$ be the angle between two unit vectors $\text{\bf v}_1,
\text{\bf v}_2\in S^2$, $0\leq\alpha\leq\pi$. Then
$$1-\cos\alpha=1-\text{\bf v}_1\cdot\text{\bf v}_2$$
is a smooth function on $S^2\times S^2$, whereas
$$\sin\alpha = |\text{\bf v}_1\times\text{\bf v}_2|$$
is not a smooth function on $S^2\times S^2$.

Similar to the excess length interpretation of $E_{\cos}$, the first 
integration in $E_{\sin}$, which is equal to
$$\int_{\gamma_\infty}\sin\alpha\,ds <\infty,$$ 
can be interpreted as the ``total momentum'' of $\gamma_\infty$ 
with respect to its asymptotic direction. Again, to minimize $E_{\sin}$, 
$\gamma$ should not have ``highly oscillated'' segments.

\bigskip

\noindent{\bf \S2. Gauss functionals}
\medskip

We define chord diagrams first, which are originated in the study of
Vassiliev knot invariants. A {\it chord diagram with n chords} consists of
$2n$ distinct points on $S^1$ which are paired into $n$ pairs. We stick a
chord to each paired points to indicate the pairing. Chord diagrams
are combinatorial objects so that two chord diagrams are thought to be the
same if they differ by an orientation preserving diffeomorphism of $S^1$
sending pairs to pairs. A {\it Gauss diagram} is a chord diagram whose end
points of chords are ordered in consistence with the orientation of $S^1$.
See Figure 2. 

\bigskip
\centerline{\epsfxsize=1in \epsfbox{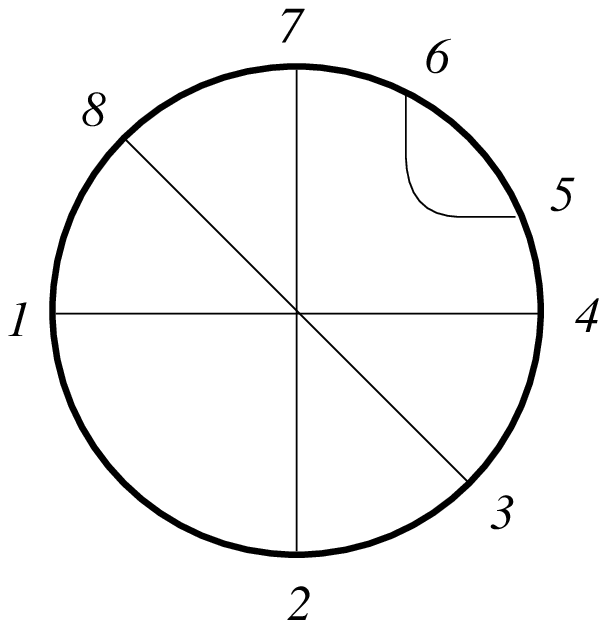}}
\smallskip
\centerline{$\ssize\text{Figure 2. A Gauss diagram.}$}
\medskip

We will denote by $C_k$ the configuration space of $k$ ordered distinct
points on $S^1$ such that the ordering of points is consistent with the
orientation of $S^1$. It is an oriented open $k$-dimensional manifold.  

Let $\omega$ be the standard volume form of the unit 2-sphere $S^2$ such that 
$$\int_{S^2}\omega=1.$$
For $x\in \Bbb R^3\setminus\{0\}$, we will denote by $\omega(x)$ the
pull-back of $\omega$ under the map
$$\Bbb R^3\setminus\{0\}\rightarrow S^2:\, x\mapsto\frac{x}{|x|}.$$
We have
$$\omega(x)=\frac1{8\pi}\,\frac{<x,dx,dx>}{|x|^3}=\frac1{4\pi}\,
\frac{x_1dx_2dx_3+x_2dx_3dx_1+x_3dx_1dx_2}{|x|^3}$$
where $<\, ,\, ,\, >$ is the mixed product in $\Bbb R^3$. 

Let $\gamma$ be a geometric knot and $D$ be a Gauss diagram with $n$
chords. We may define a functional $I_D(\gamma)$ as follows. First, we
order the chords in $D$. This allows us to define a map
$$F_D:C_{2n}\longrightarrow(S^2)^n.$$
For $(u_1,u_2,\dots,u_{2n})\in C_{2n}$, the $j$-th coordinate of its image
in $(S^2)^n$ under $F_D$ is
$$\frac{\gamma(u_{i'})-\gamma(u_{i})}{|\gamma(u_{i'})-\gamma(u_i)|}$$
if $i$ and $i'$ are paired in $D$ as the $j$-th pair and $i<i'$. Then
$$\Cal I_D(\gamma)=\int_{C_{2n}}F^*_D(\omega^n).$$
Obviously, $\Cal I_D$ is independent of the choice of orderings of chords 
in $D$ as 2-forms commute. We
will call $\Cal I_D$ a {\it Gauss functional} on geometric knots.

We may also have a {\it reduced Gauss functional} $\bar\Cal I_D$ associated 
with each chord
diagram $D$. Abuse the notation, we assume that $D$ is a Gauss diagram with 
$n$ chords. The cyclic group $\Bbb Z_{2n}$ acts on the set of Gauss diagram
having the same underlying chord diagram as that of $D$ by permuting end
points of $D$ cyclicly. Then, by definition,
$$\bar\Cal I_D=\frac1{2n}\sum_{\sigma\in\Bbb Z_{2n}}\Cal I_{\sigma D}.$$
Of course, if $\sigma D=D$ as Gauss diagrams for each $\sigma\in\Bbb Z_{2n}$,
we have $\bar\Cal I_D=\Cal I_D$.

If $D$ has only one chord, $\Cal I_D(\gamma)$ is the classical writhe integral:
$$\Cal I_w(\gamma)=\frac1{4\pi}\int_{C_{2}}\frac{<\gamma(u_2)-\gamma(u_1),
\dot\gamma(u_1),\dot\gamma(u_2)>}
{|\gamma(u_2)-\gamma(u_1)|^3}\,du_1du_2$$   
originated from Gauss' formula for the linking number of two disjoint geometric
knots in the 3-space (a link).

Suppose $D$ has two chords which cross each other. We call such a chord
diagram $X$. We have
$$\align \Cal I_X(\gamma)=\left(\frac1{4\pi}\right)^2\,\int_{C_4}
&\frac{<\gamma(u_3)-\gamma(u_1),\dot\gamma(u_1),\dot\gamma(u_3)>}
{|\gamma(u_3)-\gamma(u_1)|^3}\\
&\frac{<\gamma(u_4)-\gamma(u_2),\dot\gamma(u_2),
\dot\gamma(u_4)>}{|\gamma(u_4)-\gamma(u_2)|^3}\,du_1du_3du_2du_4.
\endalign$$

Gauss functionals contain topological information of knots.  
When geometric knots are
deformed within their isotopy classes, Gauss functionals will change their
values. But the variation of Gauss functionals may be compensated by the 
variation
of other different kinds of functionals on geometric knots so that together 
they
form {\it knot invariants}: functionals having the same value on isotopic 
geometric knots. The classical case is the so-called C\v alug\v areanu-Pohl 
self-linking formula expressing
the self-linking number of a geometric knot with nowhere vanishing curvature as
the sum of its writhe and total torsion [12]. Originated from the perturbative
expansion of Witten's Chern-Simons path integral, modern generalizations of the
self-linking formula form a vast family of knot invariants, possibly exhausts
all Vassiliev knot invariants. See [3, 5]. 

In [10], the simplest generalization of the self-linking formula was studied in
details. It involves the Gauss functional $\Cal I_X$ and another 
functional
$\Cal I_Y$. To define $\Cal I_Y$, let 
$$C_3(\gamma)=\left\{(u_1,u_2,u_3,x)\,;\,(u_1,u_2,u_3)\in C_3,\, x\in\Bbb
R^3\setminus\{\gamma(u_1),\gamma(u_2),\gamma(u_3)\}\right\}.$$ 
Also, denote
$$H(u,x)=\frac{(\gamma(u)-x)\times\dot\gamma(u)}{|\gamma(u)-x|^3}$$
For $x$ not on $\gamma$. Then
$$\Cal I_Y=\left(\frac1{4\pi}\right)^3\,\int_{C_3(\gamma)}<H(u_1,x),H(u_2,x),H(u_3,x)>\,
d^3xdu_1du_2du_3.$$
The functional
$$\frac14\Cal I_X-\frac13\Cal I_Y+\frac1{24}$$
turns out to be an integer valued knot invariant which can be identified as the
second coefficient of the Conway knot polynomial.

\bigskip
\noindent{\bf \S3. The $X$-crossing number}
\medskip

Let $D$ be a Gauss diagram with $n$ chords. For a geometric
knot $\gamma$. Define
$$\Cal C_D(\gamma)=\int_{C_{2n}}\left|F^*_D(\omega^n)\right|.$$
Here, if $F^*_D(\omega^n)=\lambda\,d_{\text{vol}}$, 
$|F^*_D(\omega^n)|=|\lambda|
\,d_{\text{vol}}$.

If $D$ is a chord diagram with $n$ chords, then define
$$\bar\Cal C_D=\frac1{2n}\sum_{\sigma\in\Bbb Z_{2n}}\Cal C_{\sigma D}.$$
We will see that $\bar\Cal C_D$ has the combinatorial meaning of being 
the ``average
$D$-crossing number''. Here, some digression about the
functionals $\bar\Cal I_w=\Cal I_w$ and $\bar\Cal C_w=\Cal C_w$ seems to be
appropriate before we could straighten out the general case.

Using a partition of unit, the 2-form $\omega$ on $S^2$ can be 
decomposed as a sum of many 2-forms supported in small
neighborhoods of single points. Let 
$\omega_{\text{\bf v}}$ be one of these 2-forms supported in a small
neighborhood of $\text{\bf v}\in S^2$ and $P_{\text{\bf v}}:\Bbb R^3\rightarrow\Bbb R^2$ be
the orthogonal projection in the direction {\bf v}. If we replace $\omega$
by $\omega_{\text{\bf v}}$   
in the functional $\Cal I_w(\gamma)$ and $\Cal C_w(\gamma)$, what we get, 
roughly speaking, are the algebraic crossing number $w(\gamma;\text{\bf v})$
and the crossing number $n(\gamma;\text{\bf v})$ of the plane projection
$P_{\text{\bf v}}(\gamma)$, respectively. To be more precise,
$$n(\gamma;\text{\bf v})=\frac12\,\sharp\{(u_1,u_2)\in C_2\,;\,
P_{\text{\bf v}}(\gamma)(u_1)=P_{\text{\bf v}}(\gamma)(u_2)\}$$
and $w(\gamma;\text{\bf v})$ represents a similar signed counting. 
The following proposition should 
therefore be
quite clear (see [7] and [6]).

\proclaim{\bf Proposition} We have
$$\Cal I_w(\gamma)=\frac1{4\pi}\,\int_{\text{\bf v}\in S^2}w(
\gamma;\text{\bf v})\,d_{\text{vol}}(\text{\bf v})$$
and
$$\Cal C_w(\gamma)=\frac1{4\pi}\,\int_{\text{\bf v}\in S^2}n(
\gamma;\text{\bf v})\,d_{\text{vol}}(\text{\bf v}),
$$
where $d_{\text{vol}}$ is the volume element of $S^2$.
\endproclaim

Therefore, $\bar\Cal C_w(\gamma)$ (so as $\bar\Cal I_w(\gamma)$) is 
the average, over all possible directions, of the crossing 
numbers that we see by looking at the geometric knot $\gamma$ 
in individual directions. Notice that we are looking at
$\gamma$ ``with one eye'' in each direction. But $\gamma$ is in the 3-space. 
To get a more stereo-scopic image, we'd better to look at it with two or more 
eyes. This is actually the principle behind stereo-photography.
With $\bar\Cal C_D$ for a general chord diagram $D$, it seems that we are
doing the same thing as in stereo-photography: First look at $\gamma$ 
through many individual eyes, and then try to combine the images obtained
individually together in a certain way. The resulting picture is of course 
more complete 
(and more complicated) than what we get by
looking at $\gamma$ with one eye.     
 
Consider now the $X$ diagram. Let $\gamma$ be a geometric knot and
$\text{\bf v}_1,\text{\bf v}_2\in S^2$. We define a number 
$n(\gamma;\text{\bf v}_1,\text{\bf v}_2)$ as follows.

We first notice that there is a subset of $S^2\times S^2$ of full measure with
the following property: If $(\text{\bf v}_1,\text{\bf v}_2)$ is in this 
subset, and $u_1,u_2,u_3,u_4\in S^1$ such that 
$$P_{\text{\bf v}_1}(\gamma)(u_1)=P_{\text{\bf v}_1}(\gamma)(u_3)\qquad
\text{and}\qquad P_{\text{\bf v}_2}(\gamma)(u_2)=P_{\text{\bf v}_2}(\gamma)
(u_4),$$
then $u_1,u_2,u_3,u_4$ are distinct. With this said, we define
$$
\align
n(\gamma;\text{\bf v}_1,\text{\bf v}_2)=\frac14
\,\sharp\{(u_1,u_2,u_3,u_4)\in C_4\,;\,
&P_{\text{\bf v}_1}
(\gamma)(u_1)=P_{\text{\bf v}_1}(\gamma)(u_3),\\
&P_{\text{\bf v}_1}(\gamma)(u_2)=P_{\text{\bf v}_1}(\gamma)(u_4)\}
\endalign
$$
for $(\text{\bf v}_1,\text{\bf v}_2)$ in that subset of full measure.

A similar argument as in the previous proposition can be used to prove 
the following proposition.

\proclaim{\bf Proposition} We have
$$\Cal C_X(\gamma)=\left(\frac1{4\pi}\right)^2\,\iint_{S^2\times S^2}
n(\gamma;\text{\bf v}_1,\text{\bf v}_2)\,d_{\text{vol}}(\text{\bf v}_1)
d_{\text{vol}}(\text{\bf v}_2).$$
\endproclaim

The {\it crossing number} of a knot is defined to be the minimum of 
crossing numbers of regular plane projections of that knot, where a plane 
projection of a knot is called {\it regular} if it has only transverse 
double pints and the number of double points is the {\it crossing number}
of that regular plane projection.

Denote by $[\gamma]$ the knot type of a geometric knot $\gamma$, and 
by $C([\gamma])$ the crossing number of the knot $[\gamma]$. Then we have
$$C([\gamma])\leq\Cal C_w(\gamma).$$
We would like to have a similar lower bound depending only on the knot type 
$[\gamma]$ for 
$\Cal C_X(\gamma)$.  

Suppose we have a plane curve with only transverse double points as its 
singular points. It is given by an immersion $S^1\rightarrow\Bbb R^2$. The 
preimages of the transverse double points are paired such that points in a
pair have the same image. This gives rise to a chord diagram. See Figure 3.

\bigskip
\centerline{\epsfxsize=1.5in \epsfbox{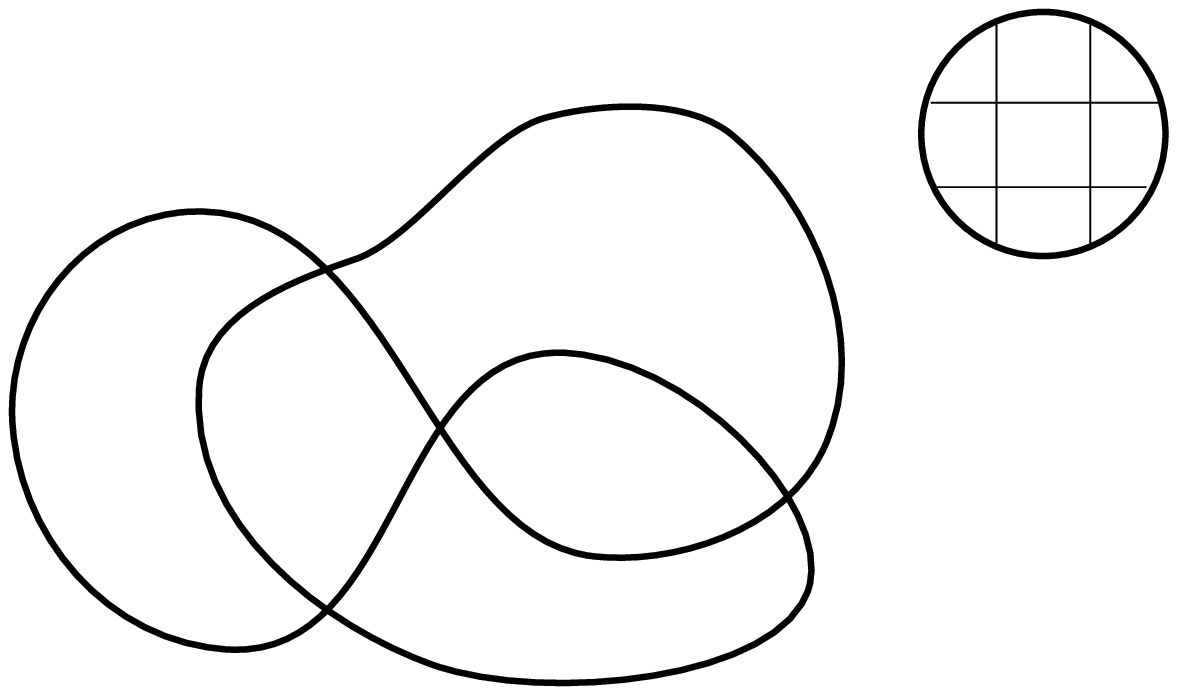}}
\smallskip
\centerline{$\ssize\text{Figure 3. A plane curve and its chord diagram.}$}
\medskip

For example, if $K$ is a knot which has the plane curve in Figure 3 as
one of its regular projection, then $C_X(K)\leq 4$.

\proclaim{\bf Definition} The {\rm $X$-crossing number} $C_X$ of a knot is the 
minimum of numbers of pairs of intersecting chords in chord diagrams of 
regular projections of that knot.
\endproclaim

\proclaim{\bf Theorem} Let $\gamma$ be a geometric knot. We have
$$C_X([\gamma])\leq\frac12\,\Cal C_X(\gamma).$$
\endproclaim  

\demo{{\bf Proof} {\rm (a sketch)}} 
The proof is based on the previous proposition and 
the fact that $\Cal C_X$ is scalar invariant. Fix a vector $\text{\bf v}_1$
such that $P_{\text{\bf v}_1}(\gamma)$ is regular. Since $\Cal C_X$ is scalar
invariant, we may assume that the preimage on $\gamma$
of the two intersecting segments of 
$P_{\text{\bf v}_1}(\gamma)$ at a double points are very close to each other. 
Then, roughly speaking, for almost all $\text{\bf v}_2$, near the double 
points of $P_{\text{\bf v}_1}(\gamma)$, we see $\gamma$ 
in the direction $\text{\bf v}_2$ as much as in the direction $\text{\bf v}_1$.
We may even see more in the direction $\text{\bf v}_2$, of course. Therefore, 
$$2\,C_X([\gamma])\leq n(\gamma;\text{\bf v}_1,\text{\bf v}_2),$$
where the factor 2 is because we don't distinguish $\text{\bf v}_1$ and
$\text{\bf v}_2$ in the counting of $C_X$. This, together with the previous 
proposition, proves the theorem. \qed
\enddemo

In general, for any chord diagram $D$, we may define the $D$-crossing number
$C_D$ of a knot and state similar theorems about $C_D$ and $\bar\Cal C_D$. We 
will not get into the details of such generalizations since they are quite
straightforward. But we do want to address one thing about 
$D$-crossing numbers unknown to the author.

A chord diagram $D$ is called a {\it sub-diagram} of another chord diagram 
$D'$ if $D$ can be obtained from $D'$ by dropping off some chords. The
{\it $D$-crossing number} of a knot is the minimum of numbers of 
$D$'s as sub-diagrams of chord diagrams of regular projections of the given 
knot. The usual crossing number is the $D$-crossing number with $D$ being
the chord diagram having one chord. It is an easy fact that the usual 
crossing number bounds the number of knots. 

\proclaim{\bf Proposition} The $X$-crossing number $C_X$ bounds the number of 
knots. In other words, given a positive number $N$, there are only finite many
knots with $C_X\leq N$.
\endproclaim

The proof is very simple. Just note that the only way to get infinitely many
chord diagrams with bounded numbers of $X$ sub-diagrams is to start with
a finite number of chord diagrams and add to them chords not intersecting with
existing ones. Chord diagrams obtained in such a way cannot be chord diagrams
of regular projections of new knots. 

\bigskip
\centerline{\epsfysize=.95in \epsfbox{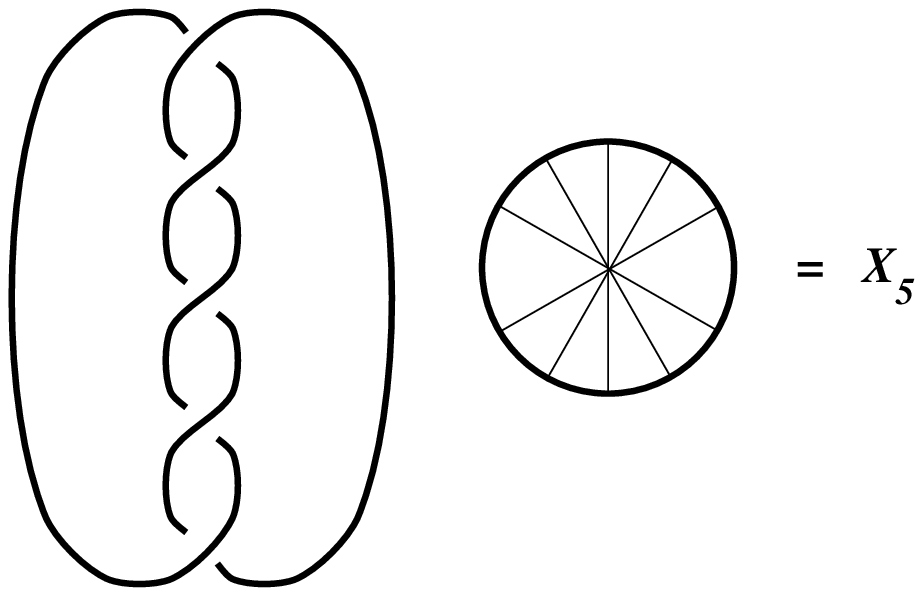}}
\smallskip
\centerline{$\ssize\text{Figure 4. A regular projection of $K_5$ 
and its chord diagram.}$}
\medskip

Let $X_n$ be the chord diagram with $n$ chords such that every pair of chords
intersect each other. For each odd $n\geq3$, we have a knot $K_n$ 
(the so-called
{\it $(2,n)$-torus knot}) which has a regular projection with $X_n$ as
its chord diagram. See Figure 4. Notice that the only 
sub-diagrams of $X_n$ are $X_k$, $k\leq
n$. Therefore, if $D$ is a chord diagram with $n$ chords and 
not equal to $X_n$, the $D$-crossing number can not bound the number 
of knots.
Moreover, as pointed out by Zhenghan Wang, the family of twist
knots in Figure 5 shows that the $X_n$-crossing number, $n\geq3$, can neither 
bound the number of knots. But it is still possible that 
$X_n$-crossing numbers could be used to control the number of knots in 
some sense. It is known that the number of knots grows at most like an 
exponential
function of the crossing number $C$. In fact, it is shown in [6] that the
number of knots with energy $E\leq M$ is at most $\simeq(0.264)(1.658)^M$. 
We say that a $D$-crossing number {\it controls} the number of knots if
the number of knots with bounded $D$-crossing numbers grows like a 
polynomial function of the crossing number. More specifically, we ask the 
following question.
                   
\proclaim{\bf Question} Does the number of knots with bounded $X_n$-crossing 
numbers grows like the function $C^{n-2}$ of the crossing number $C$? 
\endproclaim

\bigskip
\centerline{\epsfysize=1.5in \epsfbox{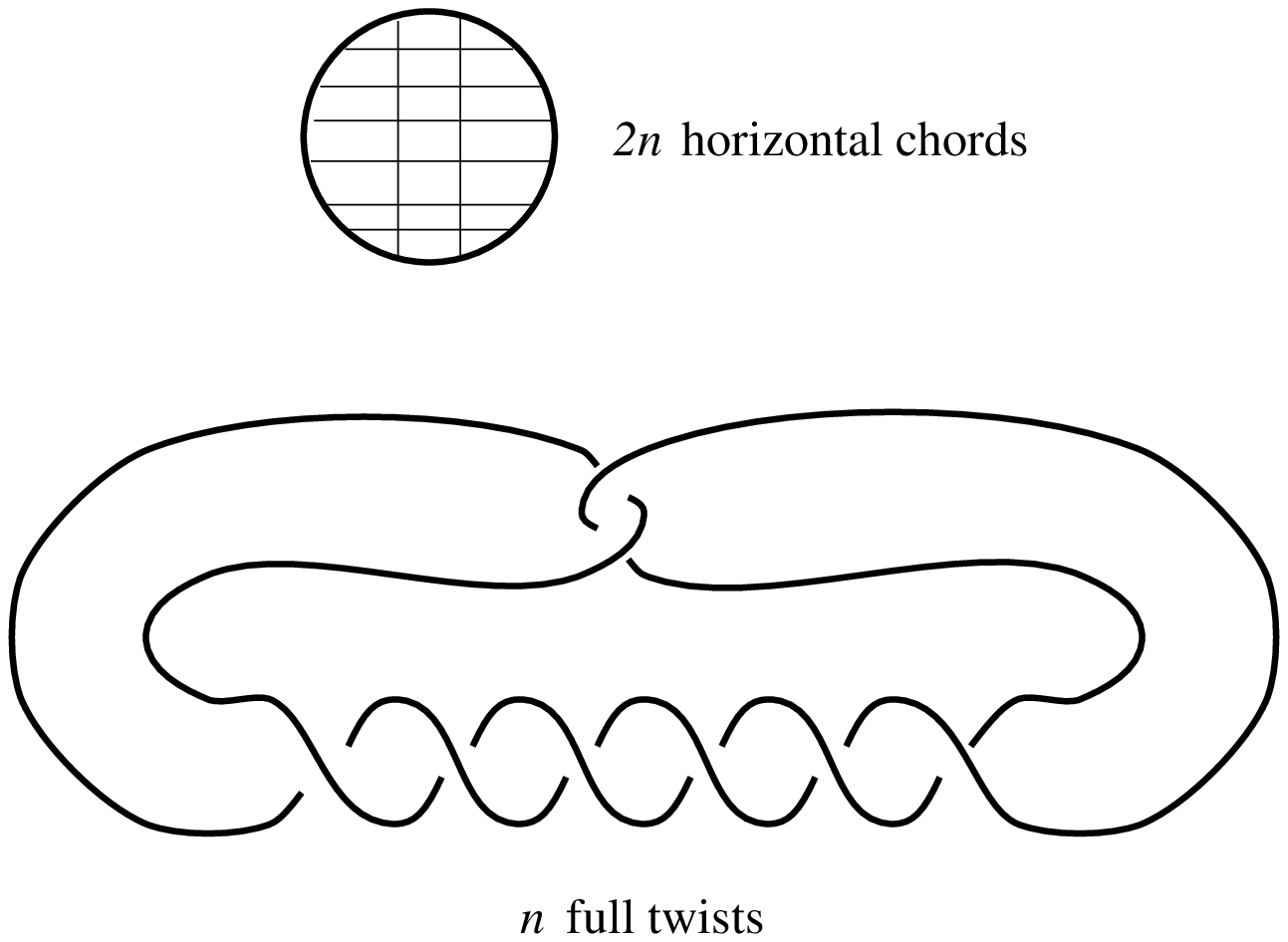}}
\smallskip
\centerline{$\ssize\text{Figure 5. Twist knots and their chord diagrams.}$}
\medskip

\bigskip
\noindent{\bf \S4. M\"obius $X$-energies}
\medskip

Here again, we will treat the chord diagram $X$ only and there is no difficult
to consider general chord diagrams.

Let $\gamma$ be a geometric knot. For $(u_1,u_2,u_3,u_4)\in C_4$, as in
\S1, let
$\alpha_{13}$ be the angle determined by $\gamma(u_1),\gamma(u_3),
\dot\gamma(u_1),\dot\gamma(u_3)$ and $\alpha_{24}$ be the angle determined by
$\gamma(u_2),\gamma(u_4),\dot\gamma(u_2),\dot\gamma(u_4)$.  We define
$$
E_{\cos,X}(\gamma)=\int_{C_4}\frac{1-\cos\alpha_{13}}
{|\gamma(u_3)-\gamma(u_1)|^2}\frac{1-\cos\alpha_{24}}
{|\gamma(u_4)-\gamma(u_2)|^2}\,d_{\text{vol}}(u_1)d_{\text{vol}}(u_2)
d_{\text{vol}}(u_3)d_{\text{vol}}(u_4)$$
and
$$
E_{\sin,X}(\gamma)=\int_{C_4}\frac{\sin\alpha_{13}}
{|\gamma(u_3)-\gamma(u_1)|^2}\frac{\sin\alpha_{24}}
{|\gamma(u_4)-\gamma(u_2)|^2}\, d_{\text{vol}}(u_1)d_{\text{vol}}(u_2)
d_{\text{vol}}(u_3)d_{\text{vol}}(u_4)$$
where $d_{\text{vol}}$ is the volume element of $\gamma$.

Following [6], we consider the following properties as essential to qualify
a functional on geometric knots being a M\"obius energy functional:
\roster
\item it is a non-negative functional and equals to zero iff on a round circle;
\item it is invariant under M\"obius transformation; and
\item it bounds the number of knots.
\endroster

\medskip
\noindent{\bf Remark:} According to the discussion of \S3, we probably
should modify (3) to ask only that general energies functionals 
{\it control} the number of knots.

The functional $E_{\cos}$ and $E_{\sin}$ obviously satisfy (1) and (2). The 
proof that (3) holds for $E_{\cos}$ is not easy. 
It depends on two main results in [6]
that $E$ is M\"obius invariant and bounds the number of knots. We will see 
below that it is much easier to prove that $E_{\sin}$ bounds the number of
knots. 

The properties (1) and (2) hold obviously for the functional $E_{\cos,X}$ and 
$E_{\sin,X}$. We have the following theorem.

\proclaim{\bf Theorem} Let $\gamma$ be a geometric knot. Then
$$C([\gamma])\leq\frac1{4\pi}\,E_{\sin}(\gamma)$$
and
$$C_X([\gamma])\leq\frac12\,\left(\frac1{4\pi}\right)^2\,E_{\sin,X}(\gamma).$$
Therefore, both $E_{\sin}$ and $E_{\sin,X}$ bound the number of knots.
\endproclaim

\demo{\bf Proof} The proof relies on some simple vector calculus. 

Suppose we have there unit vectors {\bf u, v} and {\bf w} in $\Bbb R^3$.
Consider {\bf w} as a 
segment with two vectors {\bf u} and {\bf v} stuck to its initial and
terminal points respectively. We may repeat the construction in \S1 
to get an angle $\alpha$ out of this setting, $0\leq\alpha\leq\pi$. 
If {\bf u} and {\bf w} are 
not collinear, we may let $\text{\bf u}'$ be
reflection of {\bf u} with respect to {\bf w} in the plane spanned by
{\bf u} and {\bf w}. Then $\alpha$ is the angle between $\text{\bf u}'$ 
and {\bf v}. We may calculate
$$\text{\bf u}'=2\,(\text{\bf u}\cdot\text{\bf w})\text{\bf w}-\text{\bf u}.$$
So, in particular, we have
$$<\text{\bf w},\text{\bf u},\text{\bf v}>
=-<\text{\bf w},\text{\bf u}',\text{\bf v}>.$$
Thus,
$$|<\text{\bf w},\text{\bf u},\text{\bf v}>|\leq\sin\alpha.$$
This inequality still holds if {\bf u} and {\bf w} are collinear.

For $u,v\in S^1$, use the above
inequality, we have
$$\frac{|<\gamma(v)-\gamma(u),\dot\gamma(u),\dot\gamma(v)>|}
{|\gamma(v)-\gamma(u)|^3}\leq\frac{\sin\alpha}{|\gamma(v)-\gamma(u)|^2}
|\dot\gamma(u)||\dot\gamma(v)|.$$
This implies the inequalities in the theorem. Therefor, both $E_{\sin}$ 
and $E_{\sin,X}$ bound the number of knots. \qed
\enddemo

As 
$$\sin\alpha\sim\alpha\qquad\text{and}\qquad 1-\cos\alpha\sim\frac{\alpha^2}2$$
when $\alpha$ is small, energies involving cosine appear to be ``smaller''
than those involving sine, at least locally. We don't know whether they are 
all compatible to
each other. 

\proclaim{\bf Question} Does $E_{\cos,X}$ bound the number of knots?
\endproclaim

\proclaim{\bf Question} Could we interpret $E_{\cos,X}$ as measuring 
a certain kind of excess areas?
\endproclaim

\end